\documentclass{amsart}[11pt,amssymb]
\DeclareMathSymbol{\twoheadrightarrow} {\mathrel}{AMSa}{"10}

\def\Q{{\mathbf Q}}
\def\Z{{\mathbf Z}}
\def\C{{\mathbf C}}
\def\R{{\mathbf R}}
\def\F{{\mathbf F}}

\def\R{{\mathfrak R}}

\def\CC{{\mathfrak C}}

\def\SS{{\mathbf S}}
\def\A{{\mathbf A}}
\def\Sn{{\mathbf S}_n}
\def\An{{\mathbf A}_n}

            \def\mult{\mathrm{mult}}
\def\Gal{\mathrm{Gal}}
\def\Perm{\mathrm{Perm}}

\def\Tr{\mathrm{Tr}}

\def\End{\mathrm{End}}
\def\Aut{\mathrm{Aut}}
\def\Hom{\mathrm{Hom}}

\def\cl{\mathrm{cl}}

\def\div{\mathrm{div}}
\def\Lie{\mathrm{Lie}}
\def\I{\mathrm{Id}}

\def\fchar{\mathrm{char}}

\def\GL{\mathrm{GL}}

                    \def\pr{\mathrm{pr}}

\def\M{\mathrm{M}}
\def\dim{\mathrm{dim}}
\def\O{{\mathcal O}}
\def\P{{\mathbf P}}
\def\PP{{\mathcal P}}
\def\n{{n}}

\newtheorem{thm}{Theorem}[section]
\newtheorem{lem}[thm]{Lemma}
\newtheorem{cor}[thm]{Corollary}

\theoremstyle{definition}
\newtheorem{defn}[thm]{Definition}
\newtheorem{ex}[thm]{Example}
\newtheorem{exs}[thm]{Examples}

\newtheorem{rem}[thm]{Remark}
\newtheorem{rems}[thm]{Remarks}
\title[Cyclic covers, jacobians and endomorphisms]
{Cyclic covers  of prime power degree, jacobians and endomorphisms
}
\author[Yuri G. Zarhin]{Yuri G. Zarhin}
\address{Department of Mathematics, Pennsylvania State University,
University Park, PA 16802, USA} \email{zarhin\char`\@math.psu.edu}

\begin{document}

 \begin{abstract}
Suppose $K$ is a field of characteristic zero, $K_a$ is its
algebraic closure, $f(x) \in K[x]$ is an irreducible polynomial of
degree $n \ge 5$, whose Galois group coincides either with the
full symmetric group $\Sn$ or with the alternating group $\An$.
Let $q$ be a power prime, ${\mathcal P}_q(t)=\frac{t^q-1}{t-1}$.

Let $C$ be the superelliptic curve $y^q=f(x)$ and $J(C)$ its
jacobian. We prove that if $p$ does not divide $n$ then the
algebra $\End(J(C))\otimes\Q$ of $K_a$-endomorphisms of $J(C)$ is
canonically isomorphic to $\Q[t]/{\mathcal P}_q(t)\Q[t]$.
 \end{abstract}

\maketitle

\section{Introduction}
We write $\Z,\Q,\C$  for the ring of integers, the field of
rational numbers and the field of complex numbers respectively.
Recall that a number field is called a CM-field if it is a purely
imaginary quadratic extension of a totally real field.
 Let $p$ be a prime, $q=p^r$ an integral power of p, $\zeta_q \in \C$ a primitive $q$th
root of unity, $\Q(\zeta_q)\subset \C$ the $q$th cyclotomic field
and $\Z[\zeta_q]$ the ring of integers in $\Q(\zeta_q)$. If $q=2$
then  $\Q(\zeta_q)=\Q$. It is well-known that if $q>2$ then
$\Q(\zeta_q)$ is a CM-field of degree $(p-1)p^{r-1}$. Let us put
$${\mathcal P}_q(t)=\frac{t^q-1}{t-1}=t^{q-1}+ \cdots +1 \in
\Z[t].$$ Clearly,
$$\PP(t)=\prod_{i=1}^r \Phi_{p^i}(t)$$
where
$$\Phi_{p^i}(t)=t^{(p-1)p^{i-1}}+\cdots +t^{p^{i-1}}+1\in\Z[t]$$
is the $p^i$th cyclotomic polynomial. In particular,
$$\Q[t]/\Phi_{p^i}(t)\Q[t]=\Q(\zeta_{p^i})$$ and
$$\Q[t]/{\mathcal P}_q(t)\Q[t]=\prod_{i=1}^r\Q(\zeta_{p^i}).$$

 We write
$\F_p$ for the finite field with $p$ elements.

Let $f(x) \in \C[x]$ be a polynomial of degree $n\ge 4$ without multiple roots.
 Let $C_{f,q}$ be a smooth projective model of the smooth affine curve
$$y^q=f(x).$$ Throughout this paper we assume that either $p$ does
not divide $n$ or $q$ divides $n$. It is well-known that the genus
$g(C_{f,q})$ of $C_{f,q}$ is $(q-1)(n-1)/2$ if $p$ does not divide
$n$ and
 $(q-1)(n-2)/2$ if $q$ divides $n$.
The map
 $$(x,y) \mapsto (x, \zeta_q y)$$
gives rise to a non-trivial birational automorphism
$$\delta_q: C_{f,q} \to C_{f,q}$$
of period $q$.

The jacobian $J(C_{f,q})$ of $C_{f,q}$ is an abelian variety of
dimension  $g(C_{f,q})$. We write $\End(J^{(f,q)})$ for the ring
of endomorphisms of $J^{(f,p)}$ over  $\C$ and
$\End^0(J(C_{f,q}))$ for the endomorphism algebra
$\End(J(C_{f,q}))\otimes\Q$. By Albanese functoriality, $\delta_q$
induces an automorphism of $J(C_{f,q})$ which we still denote by
$\delta_q$. One may easily check (see \ref{cycl} below) that
$$\delta_q^{q-1}+\cdots +\delta_q+1=0$$
 in $\End(J(C_{f,q}))$. This implies that if $\Q[\delta_q]$ is the
$\Q$-subalgebra of $\End^0(J(C_{f,q}))$ generated by $\delta_q$
then there is the natural surjective homomorphism
$$\Q[t]/{\mathcal P}_q(t)\Q[t] \twoheadrightarrow \Q[\delta_q]$$
which sends $t+{\mathcal P}_q(t)\Q[t]$ to $\delta_q$. One may
check that this homomorphism is, in fact, an isomorphism (see
\cite[p.~149]{Poonen}, \cite[p.~458]{SPoonen}) where the case
$q=p$ was treated).

 This gives us an embedding
$$\Q[t]/{\mathcal P}_q(t)\Q[t] \cong \Q[\delta_q] \subset \End^0(J(C_{f,q})).$$

Our main result is the following statement.

\begin{thm}
\label{endo}
Let $K$ be a subfield of $\C$ such that all the coefficients of $f(x)$
lie in $K$. Assume also that $f(x)$ is an irreducible polynomial in $K[x]$
of degree $n \ge 5$
 and its Galois group over $K$ is either the symmetric group $\Sn$ or
the alternating group $\An$. In addition, assume that either $p$
does not divide $n$ or $q\mid n$. Then
 $$\End^0(J(C_{f,q}))=\Q[\delta_q]
\cong\Q[t]/{\mathcal P}_q(t)\Q[t]=\prod_{i=1}^r\Q(\zeta_{p^i}).$$
\end{thm}

\begin{rem}
In the case when $q$ is a  prime (i.e. $q=p$) the assertion of
Theorem \ref{endo} is proven in ~\cite{ZarhinMRL,ZarhinCamb}. See
\cite{ZarhinMMJ, ZarhinNSS,{ZarhinSb}} for a discussion of finite
characteristic case.
\end{rem}

\begin{exs}
Let $n \ge 5$ be an integer, $p$ a prime, $r$ a positive integer,
$q=p^r$.
\begin{enumerate}
\item
The polynomial $x^n-x-1 \in \Q[x]$ has Galois group $\Sn$ over
$\Q$ (\cite[p.~42]{Serre}). Therefore the endomorphism algebra
(over $\C$) of the jacobian $J(C)$ of the  curve $C:y^q=x^n-x-1$
is $\Q[t]/{\mathcal P}_q(t)\Q[t]$.
\item
The  Galois group of the ``truncated exponential"
 $$\exp_n(x):=1+x+\frac{x^2}{2}+\frac{x^3}{6}+ \cdots +
\frac{x^n}{n!}\in \Q[x]$$ is either $\Sn$ or $\An$ \cite{Schur}.
Therefore the endomorphism algebra (over $\C$) of the jacobian
$J(C)$ of the  curve $C:y^q=\exp_n(x)$ is $\Q[t]/{\mathcal
P}_q(t)\Q[t]$.
\end{enumerate}
\end{exs}

\begin{rem}
\label{ribcom}
 If $f(x) \in K[x]$ then the curve $C_{f,q}$ and
its jacobian $J(C_{f,q})$ are defined over $K$. Let $K_a\subset
\C$ be the algebraic closure of $K$. Clearly, all endomorphisms of
$J(C_{f,q})$ are defined over $K_a$.
 This implies that in order to prove Theorem \ref{endo}, it suffices to
check that $\Q[\delta_q]$ coincides with  the $\Q$-algebra of
$K_a$-endomorphisms of $J(C_{f,q})$.
\end{rem}

\section{Complex abelian varieties}
\label{MT}
Throughout this section we assume that $Z$ is a complex abelian variety of positive dimension.
As usual, we write $\End^0(Z)$ for the semisimple
finite-dimensional $\Q$-algebra $\End(Z)\otimes \Q$. We write $\CC_Z$ for the center of $\End^0(Z)$.
 It is well-known that $\CC_Z$ is  a direct product of finitely many number fields.
  All the fields involved are either totally real number fields or CM-fields.
  Let $H_1(Z,\Q)$ be the first rational homology group of $Z$; it is a $2\dim(Z)$-dimensional $\Q$-vector space.
  By functoriality, $\End^0(Z)$ acts on $H_1(Z,\Q)$; hence we have an embedding
$$\End^0(Z) \hookrightarrow \End_{\Q}(H_1(Z,\Q))$$
(which sends $1$ to $1$).

Suppose $E$ is a subfield of $\End^0(Z)$ that contains the
identity map. Then $H_1(Z,\Q)$ becomes an $E$-vector space of
dimension $$d=\frac{2\dim(Z)}{[E:\Q]}.$$ We write
$$\Tr_E: \End_E(H_1(Z,\Q)) \to E$$ for the corresponding trace
map on the $E$-algebra of $E$-linear operators in $H_1(Z,\Q)$.

 Extending by
$\C$-linearity the action of $\End^0(Z)$ and of $E$ on the
complex cohomology group $$H_1(Z,\Q)\otimes_{\Q}\C=H_1(Z,\C)$$ of
$Z$ we get the embeddings $$E\otimes_{\Q}\C \subset
\End^0(Z)\otimes_{\Q}\C \hookrightarrow \End_{\C}(H_1(Z,\C))$$
which provide $H_1(Z,\C)$ with a natural structure of free
$E_{\C}:=E\otimes_{\Q}\C$-module of rank $d$.
 If $\Sigma_E$ is the set of all field embeddings  $\sigma: E \hookrightarrow \C$
 then it is well-known that
$$E_{\C}=E\otimes_{\Q}\C=\prod_{\sigma\in \Sigma_E}
E\otimes_{E,\sigma}\C=\prod_{\sigma\in \Sigma_E}\C_{\sigma}$$
where
$$\C_{\sigma}=E\otimes_{E,\sigma}\C=\C.$$
Since $H_1(Z,\C)$ is a free $E_{\C}$-module of rank $d$, there is
the corresponding trace map
$$\Tr_{E_{\C}}:\End_{E_{\C}}(H_1(Z,\C)) \to E_{\C}$$
which coincides on $E_{\C}$ with multiplication by $d$ and with
$\Tr_E$ on $\End_E(H_1(Z,\Q))$.

We write $\Lie(Z)$ for the tangent space of $Z$; it is a
$\dim(Z)$-dimensional $\C$-vector space. By functoriality,
$\End^0(Z)$ and therefore $E$ act on $\Lie(Z)$. This provides
$\Lie(Z)$ with a natural structure of $E\otimes_{\Q}\C$-module. We
have
$$\Lie(Z)=\bigoplus_{\sigma\in
\Sigma_E}\C_{\sigma}\Lie(Z)=\oplus_{\sigma\in
\Sigma_E}\Lie(Z)_{\sigma}$$ where
$$\Lie(Z)_{\sigma}=\C_{\sigma}\Lie(Z)=\{x \in \Lie(Z)\mid
ex=\sigma(e)x \quad \forall e\in E\}.$$ Let us put
$$n_{\sigma}=n_{\sigma}(Z,E)=\dim_{\C_{\sigma}}\Lie(Z)_{\sigma}=\dim_{\C}\Lie(Z)_{\sigma}.$$

\begin{rem}
\label{dual}
 Let $\Omega^1(Z)$ be the space of the differentials
of the first kind on $Z$. It is well-known that the natural map
$$\Omega^1(Z) \to \Hom_{\C}(\Lie(Z),\C)$$ is an isomorphism. This
isomorphism allows us to define via duality the natural
homomorphism
$$E  \to
\End_{\C}(\Hom_{\C}(\Lie(Z),\C))=\End_{\C}(\Omega^1(Z)).$$ This
provides $\Omega^1(Z)$ with a natural structure of
$E\otimes_{\Q}\C$-module in such a way that
$$\Omega^1(Z)_{\sigma}:=\C_{\sigma}\Omega^1(Z)\cong
\Hom_{\C}(\Lie(Z)_{\sigma},\C).$$ In particular,
$$n_{\sigma}=\dim_{\C}(\Lie(Z)_{\sigma})=\dim_{\C}(\Omega^1(Z)_{\sigma}).$$
\end{rem}

\begin{thm}
\label{mult} Suppose that $E$ contains $\CC_Z$. Then the tuple
$$(n_{\sigma})_{\sigma\in \Sigma_E} \in \prod_{\sigma\in
\Sigma_E}\C_{\sigma}=E\otimes_{\Q}\C$$ lies in
$\CC_Z\otimes_{\Q}\C$. In particular, if $E/\Q$ is Galois and
$\CC_{Z} \ne E$ then there exists a nontrivial automorphism
$\kappa: E \to E$ such that $n_{\sigma}=n_{\sigma\kappa}$ for all
$\sigma\in \Sigma_E$.
\end{thm}

\begin{proof} This is Theorem 2.3 of \cite{ZarhinCamb}.
\end{proof}

\begin{cor}
\label{cyclmult} Suppose that there exist a prime $p$, a positive
integer $r$,  the power prime  $q=p^r$ and an integer $n\ge 4$
enjoying the following properties:
\begin{itemize}
\item[(i)]
$E=\Q(\zeta_q)\subset\C$ where $\zeta_q\in\C$ is a primitive $q$th
root of unity;
\item[(ii)]
$n$ is not divisible by $p$, i.e. $n$ and $q$ are relatively
prime:
\item[(iii)]
Let $i<q$ be a positive integer that is not divisible by $p$ and
$\sigma_i:E=\Q(\zeta_q)\hookrightarrow \C$ an embedding that sends
$\zeta_q$ to $\zeta_q^{-i}$. Then
$$n_{\sigma_i}=\left[\frac{ni}{q}\right].$$
\end{itemize}

Then $\CC_{Z}=\Q(\zeta_q)$.
\end{cor}

\begin{proof}
If $q=2$ then $E=\Q(\zeta_2)=\Q$. Since $\CC_{Z}$ is a subfield of
$E=\Q$, we conclude that $\CC_{Z}=\Q=\Q(\zeta_2)$.

So, further we assume that $q>2$.
 Clearly, $\{\sigma_i\}$ is the collection $\Sigma$ of
all embeddings $\Q(\zeta_q)\hookrightarrow \C$. It is also clear
that $n_{\sigma_i}=0$ if and only if $1\le i \le [\frac{q}{n}]$.
Suppose that $\CC_{Z} \ne \Q(\zeta_q)$.
 It follows from Theorem \ref{mult}  that
there exists a non-trivial field automorphism
 $\kappa: \Q[\zeta_q] \to \Q[\zeta_q]$ such that
 for all $\sigma\in \Sigma$
$$n_{\sigma}=n_{\sigma\kappa}.$$
 Clearly, there exists an integer $m$ such that $p$ does {\sl not}
 divide $m$,
$1<m<q$ and $\kappa(\zeta_q)=\zeta_q^m$.

Assume that $q<n$. In this case the function $i \mapsto
n_{\sigma_i}=[\frac{ni}{q}]$ is strictly increasing and therefore
$n_{\sigma_i}\ne n_{\sigma_j}$ while $i\ne j$. This implies that
$\sigma_i=\sigma_i\kappa$, i.e. $\kappa$ is the identity map which
is not the case. The obtained contradiction implies that
     $$n<q.$$
Since $n \ge 4$,
$$q \ge 5.$$
 Clearly, $n_{\sigma}=0$ if and only if
$\sigma=\sigma_i$ with $1 \le i \le [\frac{q}{n}]$. Since $n$ and
$q$ are relatively prime, $[\frac{q}{n}]=[\frac{q-1}{n}]$. It
follows that $n_{\sigma_i}=0$ if and only if $1\le i \le
[\frac{q-1}{n}]$. Clearly, the map $\sigma \mapsto \sigma \kappa$
permutes the set $\{\sigma_i\mid 1 \le i \le [\frac{(q-1)}{n}], p
\mathrm{\ does\ not \ divide}\ i\}$. Since
$\kappa(\zeta_q)=\zeta_q^m$,
$\sigma_i\kappa(\zeta_q)=\zeta_q^{-im}$. This implies that
multiplication by $m$ in $(\Z/q\Z)^{*}=\Gal(\Q(\zeta_q)/\Q)$
leaves invariant the subset $$A:=\{i\bmod q\in \Z/q\Z\mid 1  \le i
\le [\frac{(q-1)}{n}], p \mathrm{\ does\ not \ divide}\ i\}.$$
Clearly, $A$ contains $1$ and therefore  $m=m\cdot 1\in A$. Since
$m<q$,
$$m=m\cdot 1\le \left[\frac{(q-1)}{n}\right] \le
\frac{(q-1)}{4}.$$ Let us consider the arithmetic progression
consisting  of  $2m$ integers $[\frac{(q-1)}{n}]+1, \ldots ,
[\frac{(q-1)}{n}]+2m$ with difference $1$. All its elements lie
between $[\frac{(q-1)}{n}]+1$ and
$$\left[\frac{(q-1)}{n}\right]+2m  \le 3\left[\frac{(q-1)}{n}\right]\le
3\frac{(q-1)}{4}<q-1.$$ Clearly, there exist exactly two elements
of  $A$ say, $d_1$ and $d_2=d_1+m$ that are  divisible by $m$.
Then  there is a positive integer $c_1$ such that
$$d_1=m c_1, d_2=m(c_1+1).$$
Clearly, either $c_1$ or $c_1+1$ is not divisible by $p$; we put
$c=c_1$ in the former case and $c=c_1+1$ in the latter case.
However, $c$ is not divisible by $p$ and
$$\left[\frac{(q-1)}{n}\right]<mc\le
\left[\frac{(q-1)}{n}\right]+2m<q-1.$$ In particular, $mc$ does
not lie in $A$. It follows that $c$ also does not lie in $A$ and
therefore
$$c>\left[\frac{(q-1)}{n}\right].$$
This means that
$$mc>m\left[\frac{(q-1)}{n}\right].$$
Since
$$mc\le\left[\frac{(q-1)}{n}\right]+2m,$$
we conclude that
$$(m-1)\left[\frac{(q-1)}{n}\right]<2m$$
and therefore
$$\left[\frac{(q-1)}{n}\right]<\frac{2m}{m-1}=2+\frac{2}{m-1}.$$
Since
$$1<m<\left[\frac{(q-1)}{n}\right],$$
we conclude that if $m>2$ then $m\ge 3$ and
$$3 \le m <\left[\frac{(q-1)}{n}\right]<2+\frac{2}{m-1}\le 3$$
and therefore $3<3$ which could not be the case. Hence $m=2$ and
$$2=m <\left[\frac{(q-1)}{n}\right]<2+\frac{2}{m-1}=4$$
and therefore
$$\left[\frac{(q-1)}{n}\right]=3.$$
It follows that $$q\ge 1+3n\ge 1+3\cdot 4=13.$$
 Since $m=2$ is not divisible by
$p$, we conclude that $p\ge 3$ and either $p=3$ and $A=\{1,2\}$ or
$p>3$ and $A=\{1,2,3\}$. In both cases $4=2\cdot 2=m\cdot 2$ must
lie in $A$. Contradiction.
\end{proof}

\section{Abelian varieties over arbitrary fields}
Let $K$ be a field. Let us fix its algebraic closure $K_a$ and
denote by $\Gal(K)$ the absolute Galois group  $\Aut(K_a/K)$ of
$K$. If $X$ is an abelian variety over $K_a$ then we write
$\End(X)$ for the ring of all its $K_a$-endomorphisms.  We write
$1_X$ (or even just $1$) for the identity automorphism of $X$. If
$Y$ is (may be another) abelian variety over $K_a$ then we write
$\Hom(X,Y)$ for the group of all $K_a$-homomorphisms from $X$ to
$Y$.
It is well-known that $\Hom(X,Y)=0$ if and only if
$\Hom(Y,X)=0$. One may easily check that if  $X$ is simple and
$\dim(X)\ge\dim(Y)$ then $\Hom(X,Y)=0$ if and only if $X$ and $Y$
are {\sl not} isogenous over $K_a$.
We write $\End^0(X)$ for the
finite-dimensional semisimple $\Q$-algebra $\End(X)\otimes\Q$ and
$\Hom^0(X,Y)$ for the finite-dimensional $\Q$-vector space
$\Hom(X,Y)\otimes\Q$. Clearly, if $X=Y$ then
$$\End^0(X)=\Hom^0(X,Y)=\Hom^0(Y,X)=\End^0(Y).$$
It is well-known that $\Hom^0(X,Y)$ and $\Hom^0(Y,X)$ have the
same dimension which does not exceed $4\dim(X)\dim(Y)$
\cite{MumfordAV}. The equality holds if and only if $\fchar(K)>0$
and both $X$ and $Y$ are supersingular abelian varieties
\cite{ZarhinMRL, ZarhinSh}.

It is well-known that if $X$ and $Y$ are  simple and the
$\Q$-algebras $\End^0(X)$ and $\End^0(Y$ are {\sl not} isomorphic
then
$$\Hom(X,Y)=0,\ \Hom(Y,X)=0.$$

 Let $E$ be a number field and $\O \subset E$ be the ring of all
its  algebraic integers. Let $(X, i)$ be a pair consisting of an
abelian variety $X$ over $K_a$ and an embedding
$$i:E \hookrightarrow  \End^0(X)$$
Here $1\in E$ must go to $1_X$. It is well known \cite{Ribet2}
that the degree $[E:\Q]$ divides $2\dim(X)$, i.e.
$$r=r_X:=\frac{2\dim(X)}{E:\Q]}$$
is a positive integer.

Let us denote by $\End^0(X,i)$ the centralizer of $i(E)$ in
$\End^0(X)$. Clearly, $i(E)$ lies in the center of the
finite-dimensional $\Q$-algebra $\End^0(X,i)$. It follows that
$\End^0(X,i)$ carries a natural structure of finite-dimensional
$E$-algebra. If $Y$ is (possibly) another abelian variety over
$K_a$ and $j:E \hookrightarrow \End^0(Y)$ is an embedding that
sends $1$ to the identity automorphism of $Y$ then we write
$$\Hom^0((X,i),(Y,j))=\{u \in \Hom^0(X,Y)\mid ui(c)=j(c)u \quad
\forall c\in E\}.$$ Clearly, $\End^0(X,i)=\Hom^0((X,i),(X,i))$. If
$d$ is a positive integer then we write $i^{(d)}$ for the
composition
$$E\hookrightarrow \End^0(X)\subset \End^0(X^d)$$
of $i$ and the diagonal inclusion $\End^0(X)\subset \End^0(X^d)$.

\begin{rem}
\label{ss}
\begin{itemize}
\item[(i)]
The $E$-algebra $\End^0(X,i)$ is semisimple. Indeed, let us split
the semisimple $\Q$-algebra $\End^0(X)$ into a finite direct
product
$$\End^0(X)= \prod_{s} D_s$$
 of simple $\Q$-algebras $D_s$. If $\pr_s:\End^0(X) \twoheadrightarrow D_s$ is the corresponding projection map and
  $D_{s,E}$ is the centralizer of $\pr_s i(E)$ in $D_s$ then one may easily check that
$$\End^0(X,i)=\prod_{s} D_{s,E}.$$
Clearly, $\pr_s i(E)\cong E$ is a simple $\Q$-algebra. It follows
from Theorem 4.3.2 on p. 104 of \cite{Herstein} that $D_{s,E}$ is
also a {\sl simple} $\Q$-algebra. This implies easily that
$D_{s,E}$ is a {\sl simple} $E$-algebra and therefore
$\End^0(X,i)$ is a semisimple $E$-algebra. It is also clear that
$\End^0(X,i)$ is a simple $E$-algebra if and only if $\End^0(X)$
is a simple $\Q$-algebra, i.e., $X$ is isogenous to a self-product
of (absolutely) simple abelian variety.
\item[(ii)]
Let $e_s$ be the identity element of $D_s$. One may view $e_s$ as
an idempotent in $\End^0(X)$. Clearly,
$$1=\sum_s e_s$$
in $\End^0(X)$ and  $e_s e_t=0$ if $s\ne t$. There exists a
positive integer $N$ such that all $N \cdot e_s$ lie in $\End(X)$.
We write $X_s$ for the image
$$X_s:=(Ne_s) (X);$$
it is an abelian subvariety in $X$ of positive dimension. Clearly,
the sum map
$$\pi_X:\prod_s X_s \to X, \quad (x_s) \mapsto \sum_s x_s$$
is an isogeny. It is also clear that the intersection $D_s\bigcap
\End(X)$ leaves $X_s \subset X$ invariant. This gives us a natural
identification
$$D_s \cong \End^0(X_s).$$
One may easily check that each $X_s$ is isogenous to a
self-product of  (absolutely) simple abelian variety. It is also
clear that
$$\Hom(X_s,X_t)=0 \quad \forall s\ne t.$$
 We write $i_s$ for the composition
$$ \pr_s i: E \hookrightarrow \End^0(X)\twoheadrightarrow
D_{s} \cong \End^0(X_s).$$ Clearly,
$$D_{s,E}=\End^0(X_s,i_s)$$ and
$$\pi_X^{-1} i \pi_X=\prod_s i_s: E \to \prod_s D_s
= \prod_s \End^0(X_s)\subset \End^0(\prod_s X_s).$$ It is also
clear that
$$\End^0(\prod_s X_s,\prod_s i_s)=\prod_sD_{s,E}.$$
\end{itemize}
\end{rem}

\begin{thm}
\label{maxE}
\begin{itemize}
\item[(i)]
$$\dim_E(\End^0((X,i)) \le
\frac{4\cdot\dim(X)^2}{[E:\Q]^2};$$
\item[(ii)]
Suppose that
 $$\dim_E(\End^0((X,i)) =
\frac{4\cdot\dim(X)^2}{[E:\Q]^2}.$$
 Then $X$ is isogenous to a self-product of
(absolutely) simple abelian variety. Also $\End^0((X,i)$ is a
central simple $E$-algebra, i.e., $E$ coincides with the center of
$\End^0((X,i)$. In addition, $X$ is an abelian variety of CM-type.

If $\fchar(K_a)=0$ then $[E:\Q]$ is even and there exist a
$\frac{[E:\Q]}{2}$-dimensional abelian variety $Z$, an isogeny
$\psi: Z^r \to X$, an embedding
$$k: E \hookrightarrow \End^0(Z)$$
that sends $1$ to  $1_Z$ and such that
$$\psi \in \Hom^0((Z^r,k^{(r)}),(X,i)).$$
\end{itemize}
\end{thm}

\begin{proof}
Recall that $r=2\dim(X)/[E:\Q]$.

 First, assume that $X$ is
isogenous to a self-product of (absolutely) simple abelian
variety, i.e., $\End^0(X,i)$ is a simple $E$-algebra. We need to
prove that
$$N:=\dim_E(\End^0(X,i))\le r^2.$$
 Let $E'$ be the center of
$\End^0(X,i)$. Let us put
$$e=[E':E].$$
Then $\End^0(X,i)$ is a {\sl central} simple $E'$-algebra of
dimension $N/e$. Then there exists a central division $E'$-algebra
$D$ such that $\End^0(X,i)$ is isomorphic to the matrix algebra
$\M_m(D)$ of size $m$ for some positive integer $m$. Dimension
arguments imply that
$$m^2\dim_{E'}(D)=\frac{N}{e},\quad \dim_{E'}(D)=\frac{N}{e m^2}.$$
Since $\dim_{E'}(D)$ is a square,
$$\frac{N}{e}=N_1^2, \quad N=e N_1^2, \quad \dim_{E'}(D)=\left(\frac{N_1}{m}\right)^2$$
for some positive integer $N_1$.
 Clearly,  $m$ divides
$N_1$.

 Clearly, $D$ contains a (maximal) field extension $L/E'$ of
degree $\frac{N_1}{m}$ and $\End^0(X,i)\cong \M_m(D)$ contains
every field extension $T/L$ of degree $m$. This implies that
$$\End^0(X) \supset\End^0(X,i)\supset T$$
and the number field $T$ has degree
$$[T:\Q]=[E':\Q]\cdot \frac{N_1}{m} \cdot m=[E:\Q]e N_1.$$
But $[T:\Q]$ must divide $2\dim(X)$; if the equality holds then
$X$ is an abelian variety of CM-type. This implies that $e N_1$
divides $r=\frac{2\dim(X)}{[E:\Q]}$. It follows that $(eN_1)^2$
divides $r^2$; if the equality holds then $X$ is an abelian
variety of CM-type. But
$$(eN_1)^2=e^2 N_1^2=e (e N_1^2)=eN=e\cdot\dim_E(\End^0(X,i)).$$
This implies that
$$\dim_E(\End^0(X,i)) \le \frac{r^2}{e} \le r^2.$$
If the equality $\dim_E(\End^0(X,i))=r^2$ holds then $e=1$ and
$$(eN_1)^2=r^2, N_1=r,\ [T:\Q]=[E:\Q]e N_1=[E:\Q] r=2\dim(X);$$
in particular, $X$ is an abelian variety of
CM-type. In addition, since $e=1$, we have $E'=E$, i.e.
$\End^0(X,i)$ is a central simple $E$-algebra.

Clearly, there exists an abelian variety $Z$ over $K_a$ with
$$E \subset D \subset \End^0(Z)$$
and an isogeny
$$\psi:Z^m \to X$$
such that the induced isomorphism  $$\End^0(Z^m) \cong \End^0(X),\
u\mapsto \psi u \psi^{-1}$$ maps identically
$$E \subset \End^0(Z)\subset \End^0(Z^m)$$
onto $E \subset \End^0(X)$.

We still have to check that if $\fchar(K)=0$ then
$$2\dim(Z)=[E:\Q].$$
 Indeed, since $D$ is a division algebra, $\dim_{\Q}(D)$ must
divide $2\dim(Z)=\frac{2\dim(X)}{m}=[E:\Q]\frac{r}{m}$. On the
other hand,
$$\dim_{\Q}(D)=[E:\Q]\dim_{E}(D)=[E:\Q]\left(\frac{r}{m}\right)^2.$$
Since  $m$ divides $r$, we conclude that $\frac{r}{m}=1$, i.e.
$$\dim_{E}(D)=1, \quad D=E, \quad 2\dim(Z)=[E:\Q].$$

Now let us consider the case of arbitrary $X$. Applying the
already proven case of the theorem to each $X_s$, we conclude that
$$\dim_E(\End^0(X,i)) \le
\left(\frac{2\dim(X_s)}{[E:\Q]}\right)^2.$$ Since
$$\End^0(X,i)=\prod_s \End^0(X_s,i_s),$$ we conclude that
$\dim_E(\End^0(X,i)) =\sum_s\dim_E(\End^0(X_s,i_s))$ does not
exceed $$\sum_s\left(\frac{2\dim(X_s)}{[E:\Q]}\right)^2\le
\frac{(2\sum_s
\dim(X_s))^2}{[E:\Q]^2}=\frac{(2\dim(X))^2}{[E:\Q]^2}.$$ It
follows that if the equality
$$\dim_E(\End^0(X,i))=\frac{(2\dim(X))^2}{[E:\Q]^2}$$ holds
then the set of indices $s$ is a singleton, i.e. $X=X_s$ is
isogenous to a self-product of  (absolutely) simple abelian
variety.
\end{proof}

Let  $d$ be a positive integer that is not divisible by
$\fchar(K)$. Let $X$ be an abelian variety
 of positive dimension
defined over $K$. We write $X_d$ for the kernel of multiplication
by $d$ in $X(K_a)$. It is known \cite{MumfordAV} that the
commutative group $X_d$ is a free $\Z/d\Z$-module of rank
$2\dim(X)$. Clearly, $X_d$ is a Galois submodule in $X(K_a)$. We
write
$$\tilde{\rho}_{d,X}:\Gal(K) \to \Aut_{\Z/d\Z}(X_d) \cong
\GL(2\dim(X),\Z/d\Z)$$ for the corresponding (continuous)
homomorphism defining the Galois action on $X_d$. Let us put
$$\tilde{G}_{d,X}=\tilde{\rho}_{d,X}(\Gal(K)) \subset
\Aut_{\Z/d\Z}(X_d).$$ Clearly,  $\tilde{G}_{d,X}$ coincides with
the Galois group of the field extension  $K(X_d)/K$ where $K(X_d)$
is the field of definition of all points on $X$ of order dividing
$d$. In particular, if a prime
 $\ell\ne \fchar(K)$  then
$X_{\ell}$ is a $2\dim(X)$-dimensional vector space over the prime
field $\F_{\ell}=\Z/\ell\Z$ and the inclusion $\tilde{G}_{\ell,X}
\subset \Aut_{\F_{\ell}}(X_{\ell})$ defines a faithful linear
representation of the group  $\tilde{G}_{\ell,X}$ in the vector
space $X_{\ell}$.

We write $\End_K(X)\subset \End(X)$ for the (sub)ring of all
$K$-endomorphisms of $X$.

Now let us assume that
$$i(\O) \subset \End_K(X).$$

 Let $\lambda$ be a maximal ideal in $\O$.
We write $k(\lambda)$ for the corresponding (finite) residue
field. Let us put
$$X_{\lambda}:=\{x \in X(K_a)\mid i(e)x=0 \quad \forall e\in \lambda\}.$$
Clearly, if $\fchar((k)(\lambda))=\ell$  then $\lambda\supset \ell
\cdot\O$ and therefore $X_{\lambda}\subset X_{\ell}$. Clearly,
$X_{\lambda}$ is a Galois submodule of $X_{\ell}$. It is also
clear that $X_{\lambda}$ carries a natural structure of
$\O/\lambda=k(\lambda)$-vector space. We write
$$\tilde{\rho}_{\lambda,X}:\Gal(K) \to
\Aut_{k(\lambda)}(X_{\lambda})$$
for the corresponding (continuous)
homomorphism defining the Galois action on $X_{\lambda}$.
 Let us
put
$$\tilde{G}_{\lambda,X}=\tilde{G}_{\lambda,i,X}:=\tilde{\rho}_{\lambda,X}(\Gal(K)) \subset
\Aut_{k(\lambda)}(X_{\lambda}).$$ Clearly, $\tilde{G}_{\lambda,X}$
coincides with the Galois group of the field extension
$K(X_{\lambda})/K$ where $K(X_{\lambda})=K(X_{\lambda,i})$ is the
field of definition of all points in $X_{\lambda}$.

In order to describe $\tilde{\rho}_{\lambda,X}$ explicitly, let us
assume for the sake of simplicity that $\lambda$ is the only
maximal ideal of $\O$ dividing $\ell$, i.e.,
$$\ell\cdot\O=\lambda^b$$
where the positive integer $b$ satisfies
$$[E:\Q]=b \cdot [k(\lambda):\F_{\ell}].$$
Then $\O\otimes\Z_{\ell}=\O_{\lambda}$ where $\O_{\lambda}$ is the
completion of $\O$ with respect to $\lambda$-adic topology. It is
well-known that that $\O_{\lambda}$ is a local principal ideal
domain and its only maximal ideal is $\lambda\O_{\lambda}$. One
may easily check that
$\ell\cdot\O_{\lambda}=(\lambda\O_{\lambda})^b$.

Let us choose an element $c \in \lambda$ that does not lie in
$\lambda^2$. Clearly, $\lambda\O_{\lambda}=c\cdot\O_{\lambda}$.
This implies that there exists a unit $u\in\O_{\lambda}^*$ such
that $\ell=u c^b$. It follows from the unique factorization of
ideals in $\O$ that
$$\lambda=\ell\cdot\O + c\cdot\O.$$
It follows readily that
$$X_{\lambda}=\{x\in X_{\ell}\mid cx=0\}\subset X_{\ell}.$$

Let $T_{\ell}(X)$ be the $\Z_{\ell}$-Tate module of $X$ defined as
projective limit of Galois modules $X_{\ell^m}$ where the
transition map(s)  $X_{\ell^{m+1}}\to X_{\ell^m}$ is
multiplication by $\ell$ \cite{MumfordAV}. Recall that
$T_{\ell}(X)$ is a free $\Z_{\ell}$-module of rank $2\dim(X)$
provided with the continuous action
$$\rho_{\ell,X}:\Gal(K) \to \Aut_{\Z_{\ell}}(T_{\ell}(X))$$
and the natural embedding
$$\End_K(X)\otimes\Z_{\ell} \hookrightarrow \End_{\Z_{\ell}}(T_{\ell}(X)),$$
whose image commutes with $\rho_{\ell,X}(\Gal(K))$. In particular,
$T_{\ell}(X)$ carries the natural structure of
$\O\otimes\Z_{\ell}=\O_{\lambda}$-module; it is known
\cite{Ribet2} that the $\O_{\lambda}$-module $T_{\ell}(X)$ is free
of rank $r=r_X=\frac{2\dim(X)}{[E:\Q]}$. There is also the natural
isomorphism of Galois modules
$$X_{\ell}=T_{\ell}(X)/\ell T_{\ell}(X),$$
which is also an isomorphism of $\End_K(X)\supset \O$-modules.
This implies that the $\O[\Gal(K)]$-module $X_{\lambda}$ coincides
with
$$c^{-1}\ell T_{\ell}(X)/\ell T_{\ell}(X)=c^{b-1}T_{\ell}(X)/c^b
T_{\ell}(X)=T_{\ell}(X)/c T_{\ell}(X)= $$ $$T_{\ell}(X)/\lambda
T_{\ell}(X) = T_{\ell}(X)/(\lambda\O_{\lambda})T_{\ell}(X).$$
Hence
$$X_{\lambda}=T_{\ell}(X)/(\lambda\O_{\lambda})T_{\ell}(X)=T_{\ell}(X)\otimes_{\O_{\lambda}}k(\lambda).$$
It follows that
$$\dim_{k(\lambda)}X_{\lambda}= \frac{2\dim(X)}{[E:\Q]}:=r_{X}.$$

Let us put
$$V_{\ell}(X)=T_{\ell}(X)\otimes_{\Z_{\ell}}\Q_{\ell};$$
it is a $2\dim(X)$-dimensional $\Q_{\ell}$-vector space that
carries a natural structure of $r_X$-dimensional
$E_{\lambda}$-vector space. There is the natural embedding
$$\End(X)\otimes\Z_{\ell}\hookrightarrow \End_{\Q_{\ell}}V_{\ell}(X).$$
Extending it by $\Q$-linearity, we get the natural embedding
$$\End^0(X)\otimes_{\Q}\Q_{\ell}\hookrightarrow
\End_{\Q_{\ell}}(V_{\ell}(X)).$$ Further we will identify
$\End^0(X)\otimes_{\Q}\Q_{\ell}$ with its image in
$\End_{\Q_{\ell}}(V_{\ell}(X))$.

\begin{rem}
\label{Edim}
 Notice that
$$E_{\lambda}=E\otimes_{\Q}\Q_{\ell}=\O\otimes\Q_{\ell}=
\O_{\lambda}\otimes_{\Z_{\ell}}\Q_{\ell}$$ is the field coinciding
with the completion of $E$ with respect to $\lambda$-adic
topology. Clearly,  $V_{\ell}(X)$ carries a natural structure of
$r_X$-dimensional $E_{\lambda}$-vector space. One may easily check
 that $\End^0(X,i)\otimes_{\Q}\Q_{\ell}$ is a
$E\otimes_{\Q}\Q_{\ell}=E_{\lambda}$-vector subspace (even
subalgebra) in  $\End_{E_{\lambda}}(V_{\ell}(X))$. Clearly,
$$\dim_{E_{\lambda}}(\End^0(X,i)\otimes_{\Q}\Q_{\ell})=\dim_E(\End^0(X,i))$$
and
$$\dim_{E_{\lambda}}(\End_{E_{\lambda}}(V_{\ell}(X)))=r_X^2.$$
This implies that
$$\End^0(X,i)\otimes_{\Q}\Q_{\ell}=\End_{E_{\lambda}}(V_{\ell}(X))$$
if and only if
$$\dim_E(\End^0(X,i))=r_X^2.$$
\end{rem}

 Using the inclusion
$$\Aut_{\Z_{\ell}}(T_{\ell}(X)) \subset
\Aut_{\Q_{\ell}}(V_{\ell}(X)),$$ one may view $\rho_{\ell,X}$ as
$\ell$-adic representation
$$\rho_{\ell,X}:\Gal(K) \to \Aut_{\Z_{\ell}}(T_{\ell}(X))\subset
\Aut_{\Q_{\ell}}(V_{\ell}(X)).$$

Since $X$ is defined over $K$, one may associate with every $u \in
\End(X)$ and $ \sigma \in \Gal(K)$ an endomorphism $^{\sigma}u \in
\End(X)$ such that $$^{\sigma}u(x)=\sigma u(\sigma^{-1}x) \quad
\forall x \in X(K_a).$$ Clearly,
$$^{\sigma}u=u \quad \forall u \in \End_K(X).$$
In particular,
$$^{\sigma}u=u \quad \forall u \in \O$$
(here we identify $\O$ with $i(\O)$). It follows easily that for
each $ \sigma \in \Gal(K)$  the map $u \to ^{\sigma}u$ extends by
$\Q$-linearity to a certain automorphism of $\End^0(X)$. It is
also clear that $^{\sigma}u =u$ for each $u\in E$ and
$$^{\sigma}u \in \End^0(X,i) \quad \forall u \in \End^0(X,i),\sigma \in
\Gal(K).$$

\begin{rem}
\label{tate}
 The definition of $T_{\ell}(X)$ as  the projective
limit of Galois modules $X_{\ell^m}$  implies that
$$^{\sigma}u(x)=\rho_{\ell,X}(\sigma) u \rho_{\ell,X}(\sigma)^{-1}(x) \quad
\forall x \in T_{\ell}(X).$$ It follows easily that
$$^{\sigma}u(x)=\rho_{\ell,X}(\sigma) u \rho_{\ell,X}(\sigma)^{-1}(x) \quad
\forall x \in V_{\ell}(X), u \in \End^0(X),\sigma \in \Gal(K).$$
This implies that for each $\sigma \in \Gal(K)$
$$\rho_{\ell,X}(\sigma) \in \Aut_{E_{\lambda}}(V_{\lambda}(X)).$$
and therefore
$$\rho_{\ell,X}(\Gal(K))\subset \Aut_{E_{\lambda}}(V_{\lambda}(X))$$
\cite{SerreAb,Ribet2}. It is also clear that
$$\rho_{\ell,X}(\sigma) u \rho_{\ell,X}(\sigma)^{-1}\in \End^0(X)\otimes_{\Q}\Q_{\ell}
\quad \forall u\in \End^0(X)\otimes_{\Q}\Q_{\ell}$$ and
$$\rho_{\ell,X}(\sigma) u \rho_{\ell,X}(\sigma)^{-1}\in \End^0(X,i)\otimes_{\Q}\Q_{\ell}
\quad \forall u\in \End^0(X,i)\otimes_{\Q}\Q_{\ell}.$$
\end{rem}

We refer to ~\cite{ZarhinTexel,ZarhinMRL2,ZarhinMMJ,ZarhinVery}
for a discussion of the following definition.

\begin{defn}
Let $V$ be a vector space over a field $\F$, let $G$ be a group
and $\rho: G \to \Aut_{\F}(V)$ a linear representation of $G$ in
$V$. We say that the $G$-module $V$ is {\sl very simple} if it
enjoys the following property:

If $R \subset \End_{\F}(V)$ is an $\F$-subalgebra containing the
identity operator $\I$ such that

 $$\rho(\sigma) R \rho(\sigma)^{-1} \subset R \quad \forall \sigma \in G$$
 then either $R=\F\cdot \I$ or $R=\End_{\F}(V)$.
\end{defn}

\begin{rems}
\label{image}
\begin{enumerate}
\item[(i)]
If $G'$ is a subgroup of $G$ and the $G'$-module $V$ is very
simple then obviously the $G$-module $V$ is also very simple.

\item[(ii)]
Clearly, the $G$-module $V$ is very simple if and only if the
corresponding $\rho(G)$-module $V$ is very simple. This implies
easily that if $H \twoheadrightarrow G$ is a surjective group
homomorphism then the $G$-module $V$ is  very simple if and only
if the corresponding $H$-module $V$ is  very simple.

\item[(iii)]
Let $G'$ be a normal subgroup of $G$. If $V$ is a very simple
$G$-module then either $\rho(G') \subset \Aut_{k}(V)$ consists of
scalars (i.e., lies in $k\cdot\I$) or the $G'$-module $V$ is
absolutely simple. See ~\cite[Remark 5.2(iv)]{ZarhinMMJ}.

\item[(iv)]
Suppose $F$ is a discrete valuation field with  valuation ring
$O_F$, maximal ideal $m_F$ and residue field $k=O_F/m_F$. Suppose
$V_F$ a finite-dimensional $F$-vector space, $\rho_F: G \to
\Aut_F(V_F)$ a $F$-linear representation of $G$. Suppose $T$ is a
$G$-stable $O_F$-lattice in $V_F$ and the corresponding
$k[G]$-module $T/m_F T$ is isomorphic to $V$. Assume that the
$G$-module $V$ is very simple. Then the $G$-module $V_F$ is also
very simple. See ~\cite[Remark 5.2(v)]{ZarhinMMJ}.
\end{enumerate}
\end{rems}

\begin{thm}
\label{Very}
Suppose that $X$ is an abelian variety defined over $K$ and
$i(\O)\subset \End_K(X)$. Let $\ell$ be a prime different from
$\fchar(K)$. Suppose that $\lambda$ is the only maximal ideal
dividing $\ell$ in $\O$. Suppose that the natural representation
in the $k(\lambda)$-vector space $X_{\lambda}$ is very simple.
Then $\End^0(X,i)$ enjoy one of the following two properties:
\begin{itemize}
\item[(i)]
$\End^0(X,i)=i(E)$, i.e. $i(E)\cong E$ is a maximal commutative
subalgebra in $\End^0(X)$ and $i(\O)\cong \O$ is a maximal
commutative subring in $\End(X)$;
\item[(ii)]
$\End^0(X,i)$ is a central simple $E$-algebra of dimension $r_X^2$
and $X$ is an abelian variety of CM-type over $K_a$. In addition,
if $\fchar(K)=0$ then $[E:\Q]$ is even and there exist a
$\frac{[E:\Q]}{2}$-dimensional abelian variety $Z$, an isogeny
$\psi: Z^r \to X$ and an embedding
$$k: E \hookrightarrow \End^0(Z)$$
that sends $1$ to  $1_Z$  such that
$$\psi \in \Hom^0((Z^r,k^{(r)}),(X,i)).$$
\end{itemize}
\end{thm}

\begin{proof}
In light of \ref{image}(ii), the $\Gal(K)$-module $X_{\lambda}$ is
very simple. In light of \ref{image}(iv) and Remark \ref{tate}
$$\rho_{\ell,X}:\Gal(K) \to \Aut_{E_{\lambda}}(V_{\ell}(X))$$
is also very simple.
Let us put
$$R=\End^0(X,i)\otimes_{\Q}\Q_{\ell}.$$
It follows from Remark \ref{tate} that either $R=E_{\lambda}\I$ or
$R=\End_{E_{\lambda}}(V_{\ell}(X))$. By Remark \ref{Edim},
$$\dim_{E_{\lambda}}(R)=\dim_{E_{\lambda}}(\End^0(X,i)\otimes_{\Q}\Q_{\ell})=
\dim_E(\End^0(X,i)).$$ It follows that  $\dim_E(\End^0(X,i))=1$ or
$r_X^2$. Clearly, if $\dim_E(\End^0(X,i))=1$ then
$\End^0(X,i)=i(E)$ and the property (i) holds. Suppose that
$\dim_E(\End^0(X,i))=r_X^2$. Applying Theorem \ref {maxE}, we
conclude that the property (ii) holds.
\end{proof}

Let $Y$ be an abelian variety of positive dimension over $K_a$ and
$u$ a non-zero endomorphism of $Y$. Let us consider the abelian
(sub)variety
$$Z=u(Y) \subset Y.$$

\begin{rem} \label{rat} If $Y$ is defined over $K$ and
$u\in\End_K(Y)$ then $Z$ is also defined over $K$ and the
inclusion map $Z\subset Y$ is defined over $K$. Indeed, clearly,
$Z$ and the inclusion map $Z\subset Y$ are defined over
$K_a^{\Gal(K)}$, i.e. $Z$ and $Z\subset Y$ are defined over a
purely inseparable extension of $K$. By Theorem of Chow \cite[Th.
5 on p. 26]{Lang}, $Z$ is defined over $K$. It follows that every
homomorphism between $Z$ and $Y$ is defined over a separable
extension of $K$. Hence $Z\subset Y$ is defined over $K$.
\end{rem}

We write $\Omega^1(Y)$ (resp. $\Omega^1(Z)$) for the
$\dim(Y)$-dimensional (resp. $\dim(Z)$-dimensional) $K_a$-vector
space of differentials of the first kind on $Y$ (resp. on $Z$).

\begin{thm}
\label{chardif} Let $Y$ be an abelian variety of positive
dimension over $K_a$ and $\delta$  an automorphism of $Y$. Suppose
that the induced $K_a$-linear operator
$$\delta^*:\Omega^1(Y) \to \Omega^1(Y)$$
is diagonalizable. Let  $S$ be the set of eigenvalues of $\delta^*$ and
$\mult_Y:S \to \Z_+$  the integer-valued function which assigns
to each eigenvalue its multiplicity.

Suppose that $P(t)$ is a polynomial with integer coefficients such
that $u=P(\delta)$ is a non-zero endomorphism of $Y$. Let us put
$Z=u(Y)$. Clearly, $Z$ is $\delta$-invariant and we write
$\delta_Z:Z \to Z$ for the corresponding automorphism of $Z$ (i.e. for the restriction of $\delta$ to $z$).
Suppose that
$$\dim(Z)=\sum_{\lambda \in S, P(\lambda)\ne 0}\mult_Y(\lambda).$$
Then the spectrum of $\delta_Z^*:\Omega^1(Z) \to \Omega^1(Z)$
coincides with $S_P=\{\lambda \in S, P(\lambda)\ne 0\}$ and the
multiplicity of an eigenvalue $\lambda$ of $\delta_Z^*$ equals
$mult_Y(\lambda)$.
\end{thm}

\begin{proof}
Clearly, $u$ commutes with $\delta$. We write $v$ for the
(surjective) homomorphism $Y \twoheadrightarrow Z$ induced by $u$
and $j$ for the inclusion map $Z\subset Y$.  Notice that $u:Y\to
Y$ splits into a composition
$$Y \stackrel{v}{\twoheadrightarrow} Z \stackrel{j}{\hookrightarrow}
Y,$$ i.e. $u=jv$. Clearly,
$$\delta_Z v=v\delta \in \Hom(Y,Z),\quad j\delta_Z=\delta j\in
\Hom(Z,Y),\quad u=jv\in\End(Y), \quad u\delta=\delta u\in
\End(Y).$$ It is also clear that the induced map
$$u^*:\Omega^1(Y) \to \Omega^1(Y)$$
coincides with $P(\delta^*)$. It follows that
$$u^*(\Omega^1(Y))=P(\delta^*)(\Omega^1(Y))$$
has dimension
$$\sum_{\lambda \in S, P(\lambda)\ne 0}\mult_Y(\lambda)=\dim(Y)$$
and coincides with $$\oplus_{\lambda \in S, P(\lambda)\ne
0}W_{\lambda}$$ where $W_{\lambda}$ is the eigenspace of $\delta$
attached to eigenvalue $\lambda$.
 Since $u^*=v^*j^*$,
$$u^*(\Omega^1(Y))=v^*j^*(\Omega^1(Y))\subset v^*(\Omega^1(Z)).$$
Since
$$\dim(u^*(\Omega^1(Y)))=\dim(Y)=\dim(\Omega^1(Z))\ge
\dim(v^*(\Omega^1(Z))),$$ the subspace
$$u^*(\Omega^1(Y))=v^*(\Omega^1(Z))$$ and
$$v^*:\Omega^1(Z)\hookrightarrow \Omega^1(Y).$$ It follows that if
we denote by $w$ the isomorphism $v^*:\Omega^1(Z)\cong
v^*(\Omega^1(Z))$  and by $\gamma$ the restriction of
$\delta^*$ to $v^*(\Omega^1(Z))$ then $\gamma w=w \delta_Y^* $ and
therefore
$$\gamma=w \delta_Y^* w^{-1}.$$
\end{proof}

\section{Cyclic covers and jacobians}
\label{Prelim}

 Throughout this paper we fix a prime number $p$ and its integral power $q=p^r$ and assume
that $K$ is a field of characteristic different from $p$. We fix
an algebraic closure $K_a$ and write $\Gal(K)$ for the absolute
Galois group $\Aut(K_a/K)$. We also fix in $K_a$ a primitive $q$th
root of unity $\zeta$.

Let $f(x) \in K[x]$ be a separable polynomial of degree $n \ge 4$.
We write $\R_f$ for the set of its roots and denote by
$L=L_f=K(\R_f)\subset K_a$ the corresponding splitting field. As
usual, the Galois group $\Gal(L/K)$ is called the Galois group of
$f$ and denoted by $\Gal(f)$. Clearly, $\Gal(f)$ permutes elements
of $\R_f$ and the natural map of $\Gal(f)$ into the group
$\Perm(\R_f)$ of all permutations of $\R_f$ is an embedding. We
will identify $\Gal(f)$ with its image and consider it as a
permutation group of $\R_f$. Clearly, $\Gal(f)$ is transitive if
and only if $f$ is irreducible in $K[x]$.

Further, we assume that either $p$ does not divide $n$ or $q$ does
divide $n$.

If $p$ does not divide $n$ then we write (as in
\cite{ZarhinCrelle})
$$V_{f,p}=(\F_p^{\R_f})^{00}=(\F_p^{\R_f})^{0}$$ for the $(n-1)$-dimensional
$\F_p$-vector space of functions
$$\phi:\R_f \to \F_p, \ \sum_{\alpha\in \R_f}\phi(\alpha)=0\}$$
provided with a natural action of the permutation group
$\Gal(f)\subset \Perm(\R_f)$. It is  the {\sl heart} over the
field $\F_p$ of the group $\Gal(f)$ acting on the set $\R_f$
 \cite{Mortimer,ZarhinCrelle}.

 \begin{rem}
\label{SnAn} If $p$ does not divide $n$ and $\Gal(f)=\Sn$ or $\An$
then the $\Gal(f)$-module $V_{f,p}$ is very simple.
 \end{rem}

 Let $C=C_{f,q}$ be the smooth projective model of the smooth
affine $K$-curve
            $$y^q=f(x).$$
So $C$ is a smooth projective curve defined over $K$. The
rational function $x \in K(C)$ defines a finite cover  $\pi:C \to
\P^1$ of degree $p$. Let $B'\subset C(K_a)$ be the set of
ramification points.  Clearly, the restriction of $\pi$ to $B'$ is
an {\sl injective} map $B' \hookrightarrow \P^1(K_a)$, whose image
is the disjoint union of $\infty$ and  $\R_f$ if $p$ does {\sl
not} divide $\deg(f)$ and just $\R_f$ if it does. We write
$$B=\pi^{-1}(\R_f)=\{(\alpha,0)\mid \alpha \in \R_f\} \subset B' \subset C(K_a).$$
Clearly, $\pi$ is ramified at each point of $B$ with ramification
index $q$. We have $B'=B$ if and only if $n$ is  divisible by $p$.
If $n$ is not divisible by $p$ then $B'$ is the disjoint union of
$B$ and a single point $\infty':=\pi^{-1}(\infty)$. In addition,
the ramification index of $\pi$ at $\pi^{-1}(\infty)$ is also $q$.
Using Hurwitz's formula, one may easily compute the genus
$g=g(C)=g(C_{q,f})$ of $C$ (\cite[pp.~401--402]{Koo},
\cite[proposition 1 on p. 3359]{Towse}, \cite[p. 148]{Poonen}).
Namely, $g$ is $(q-1)(n-1)/2$ if $p$ does {\sl not} divide $n$ and
$(q-1)(n-2)/2$ if $q$ does divide $n$.

\begin{rem}
\label{genre}
 Assume that $p$ does not divide $n$ and consider  the plane triangle
(Newton polygon)
$$\Delta_{n,q}:=\{(j,i)\mid 0\le j,\quad 0\le i, \quad qj+ni\le nq\}$$
with the vertices $(0,0)$, $(0,q)$ and $(n,0)$. Let $L_{n,q}$ be
the set of integer points in the interior of $\Delta_{n,q}$.
 One may easily check that $g=(q-1)(n-1)/2$ coincides with the number
of elements of $L_{n,q}$.  It is also clear that for each
$(j,i)\in L_{n,q}$
$$1\le j \le n-1; \quad 1 \le i \le q-1;\quad q(j-1)+(j+1)\le n(q-i).$$
Elementary calculations (\cite[theorem 3 on p. 403]{Koo}) show
that
$$\omega_{j,i}:=x^{j-1}dx/y^{q-i}=x^{j-1}y^idx/y^q=x^{j-1}y^{i-1} dx/y^{q-1}$$ is a differential of the first
kind on $C$ for each $(j,i) \in L_{n,q}$. This implies easily that
the collection $\{\omega_{j,i}\}_{(j,i)\in L_{n,q}}$ is a basis in
the space of differentials of the first kind on $C$.
\end{rem}

 There is a non-trivial birational
 $K_a$-automorphism of $C$
 $$\delta_q:(x,y) \mapsto (x, \zeta y).$$
Clearly, $\delta_q^q$ is the identity map and
 the set of fixed points of $\delta_q$ coincides with $B'$.

 \begin{rem}
\label{nondiv} Let us assume that $n=\deg(f)$ is divisible by $q$
say, $n=qm$ for some positive integer $m$.
 Let $\alpha \in K_a$ be a root of $f$ and $K_1=K(\alpha)$ be
the corresponding subfield of $K_a$. We have
$f(x)=(x-\alpha)f_1(x)$ with $f_1(x) \in K_1[x]$. Clearly,
$f_1(x)$ is a separable polynomial over $K_1$ of   degree
$qm-1=n-1 \ge 4$. It is also clear that the polynomials
$$h(x)=f_1(x+\alpha), h_1(x)=x^{n-1}h(1/x) \in K_1[x]$$
are separable of the same degree $qm-1=n-1\ge 4$. The standard
substitution $$x_1=1/(x-\alpha), y_1=y/(x-\alpha)^m$$ establishes
a birational isomorphism between $C_{f,p}$ and a curve
$$C_{h_1}: y_1^q=h_1(x_1)$$ (see \cite[p.~3359]{Towse}). In
particular, the jacobians of $C_f$ and $C_{h_1}$ are isomorphic
over $K_a$  (and even over $K_1$). But $\deg(h_1)=qm-1$ is {\sl
not} divisible by $p$. Clearly, this isomorphism commutes with the
actions of $\delta_q$. Notice also that if the Galois group of $f$
over $K$ is $\Sn$ (resp. $\An$) then the Galois group of $h_1$
over $K_1$ is $\SS_{n-1}$ (resp. $\A_{n-1}$).
\end{rem}

\begin{rem}
\begin{itemize}
\item[(i)]
\label{dfk} Let $\Omega^1(C)=\Omega^1(C_{(f,q)})$ be the
$K$-vector space of differentials of the first kind on $C$. It is
well-known that $\dim_K(\Omega^1(C_{(f,q)}))$ coincides with the
genus of $C_{(f,q)}$. By functoriality, $\delta_q$ induces on
$\Omega^1(C_{(f,q)})$  a certain $K$-linear automorphism
$$\delta_q^*:\Omega^1(C_{(f,q)})\to \Omega^1(C_{(f,q)}).$$
 Clearly, if for some positive integer $j$ the
differential $\omega_{j,i}=x^{j-1} dx/y^{q-i}$ lies in
$\Omega^1(C_{(f,q)})$ then it is an eigenvector of $\delta_q^*$
with eigenvalue $\zeta^{i}$.
\item[(ii)]
Now assume that $p$ does {\sl not} divide $n$. It follows from
Remark \ref{genre} that the collection
$$\{\omega_{j,i}=x^{j-1} dx/y^{q-i}\mid (i,j) \in L_{n,q}\}$$
is an eigenbasis of $\Omega^1(C_{(f,q)})$. This implies that the
multiplicity of the eigenvalue $\zeta^{-i}$ of $\delta_q^*$
coincides with number of interior integer points in $\Delta_{n,q}$
along the corresponding (to $q-i$) horizontal line. Elementary
calculations show that this number is $\left[\frac{ni}{q}\right]$;
in particular, $\zeta^{-i}$ {\sl is an eigenvalue  if and only if}
$\left[\frac{ni}{q}\right]>0$. Taking into account that $n \ge 4$
and $q=p^r$, we conclude that $\zeta^i$ is an eigenvalue of
$\delta_q^*$ for each integer $i$ with $p^r-p^{r-1}\le i \le
p^r-1=q-1$. It also follows easily that $1$ is {\sl not} an
eigenvalue $\delta_q^*$. This implies that
$${\mathcal P}_q(\delta_q^*)={\delta_q^*}^{q-1}+\cdots +\delta_q^*+1=0$$
in $\End_K(\Omega^1(C_{(f,q)}))$. In addition, one may easily
check that if ${\mathcal H}(t)$ is a polynomial with rational
coefficients such that ${\mathcal H}(\delta_q^*)=0$ in
$\End_K(\Omega^1(C_{(f,q)}))$ then ${\mathcal H}(t)$ is {\sl
divisible} by ${\mathcal P}_q(t)$ in $\Q[t]$. \end{itemize}
\end{rem}

Let $J(C_{f,q})=J(C)=J(C_{f,q})$ be the jacobian of $C$. It is a
$g$-dimensional abelian variety defined
 over $K$ and one may view (via Albanese functoriality) $\delta_q$ as an element of
 $$\Aut(C) \subset\Aut(J(C)) \subset \End(J(C))$$
such that
  $\delta_q \ne \I$ but $\delta_q^q=\I$
where $\I$ is the identity endomorphism of $J(C)$. Here $\Aut(C)$
stands for the group of $K_a$-automorphisms of $C$, $\Aut(J(C))$
stands for the group of $K_a$-automorphisms of $J(C)$ and
$\End(J(C))$ stands for the ring of all $K_a$-endomorphisms of
$J(C)$. We write $\Z[\delta_q]$ for the subring of $\End(J(C))$
generated by $\delta_q$.
 As usual, we write $\End^0(J(C))=\End^0(J(C_{f,q}))$ for
the corresponding $\Q$-algebra $\End(J(C))\otimes \Q$. We write
$\Q[\delta_q]$ for the $\Q$-subalgebra of $\End^0(J(C))$ generated
by $\delta_q$.

\begin{rem}
\label{alb} Assume that $p$ does not divide $n$.
 Let $P_0$ be one of the $\delta_q$-invariant points
(i.e., a ramification point for $\pi$) of $C_{f,p}(K_a)$. Then
$$\tau: C_{f,q} \to J(C_{f,q}), \quad P\mapsto \cl((P)-(P_0))$$ is
an embedding of complex algebraic varieties and it is well-known
that the induced map $$\tau^*: \Omega^1(J(C_{f,q})) \to
\Omega^1(C_{f,q})$$ is a $\C$-linear isomorphism obviously
commuting with the actions of $\delta_q$. (Here $\cl$ stands for
the linear equivalence class.) This implies that $n_{\sigma_i}$
coincides with the dimension of the eigenspace of
$\Omega^1(C_{(f,q)})$ attached to the eigenvalue $\zeta^{-i}$ of
$\delta_q^*$. Applying Remark \ref{dfk}, we conclude that if
${\mathcal H}(t)$ is a monic polynomial with integer coefficients
such that ${\mathcal H}(\delta_q)=0$ in $\End(J^{(f,q)})$ then
${\mathcal H}(t)$ is divisible by ${\mathcal P}_q(t)$ in $\Q[t]$
and therefore in $\Z[t]$.
\end{rem}

\begin{rem}
\label{multprim} Assume that $p$ does not divide $n$.  Clearly,
the set $S$ of eigenvalues $\lambda$ of
$$\delta_q^*:\Omega^1(J(C_{f,q})) \to \Omega^1(J(C_{f,q}))$$
with $\PP_{q/p}(\lambda)\ne 0$  consists of {\sl primitive} $q$th
roots of unity $\zeta^{-i}$ ($1 \le i<q, (i,p)=1$) with
$\left[\frac{ni}{q}\right]>0$ and the multiplicity of $\zeta^{-i}$
equals  $\left[\frac{ni}{q}\right]$, thanks to Remarks \ref{alb}
and \ref{dfk}. Let us compute the sum
$$M=\sum_{1 \le i<q,(i,p)=1}\left[\frac{ni}{q}\right]$$
 of multiplicities of eigenvalues from $S$.

First, assume that  $q>2$. Then $\varphi(q)=(p-1)p^{r-1}$ is even
and  for each  (index) $i$ the difference $q-i$ is also prime to
$p$, lies between $1$ and $q$ and
$$\left[\frac{ni}{q}\right]+\left[\frac{n(q-i)}{q}\right]=n-1.$$
It follows easily that
$$M=(n-1)\frac{\varphi(q)}{2}=\frac{(n-1)(p-1)p^{r-1}}{2}.$$
Now assume that $q=p=2$ and therefore $r=1$. Then $n$ is odd,
$$C_{f,q}=C_{f,2}: y^2=f(x)$$
 is a hyperelliptic curve of genus $g=\frac{n-1}{2}$ and
$$\delta_2:(x,y) \mapsto (x,-y).$$
It is well-known that  the differentials $x^i \frac{dx}{y}$ $(0\le
i\le g-1$) constitute a basis of the $g$-dimensional
$\Omega^1(J(C_{f,2}))$. It follows that $\delta_2^*$ is just
multiplication by $-1$.  Therefore
$$M=g=\frac{n-1}{2}=\frac{(n-1)(p-1)p^{r-1}}{2}.$$

Notice  that if the abelian (sub)variety
$Z:=\PP_{q/p}(\delta_q)(J(C_{f,q}))$ has dimension $M$ then the
data $Y=J(C_{f,q}), \delta=\delta_q, P=\PP_{q/p}(t)$ satisfy the
conditions of Theorem \ref{chardif}.
\end{rem}

\begin{lem}
\label{order} Assume that $p$ does not divide $n$.  Let
$D=\sum_{P\in B}a_P (P)$ be a divisor on $C=C_{f,p}$ with degree
$0$ and support in $B$. Then $D$ is principal if and only if all
the coefficients $a_P$ are divisible by $q$.
\end{lem}

\begin{proof}
Suppose $D=\div(h)$ where $h\in K_a(C)$ is a non-zero rational
function of $C$. Since $D$ is $\delta_q$-invariant, the rational
function
$$\delta_q^* h:=h\delta_q= c \cdot h$$
for some non-zero $c \in K_a$. It follows easily from the
$\delta_q$-invariance of the splitting
$$K_a(C)=\oplus_{i=0}^{q-1}y^i\cdot K_a(x)$$
that
$$h=y^i \cdot u(x)$$
for some non-zero rational function $u(x) \in K_a(x)$ and a
non-negative integer $i\le q-1$. It follows easily that all finite
zeros and poles of $u(x)$ lie in $B$, i.e., there exists an
integer-valued function $b$ on $\R_f$ such that $u$ coincides, up
to multiplication by a non-zero constant, to $\prod_{\alpha\in
\R_f}(x-\alpha)^{b(\alpha)}$. Notice that
$$\div(y)=\sum_{P\in B}(P)- n(\infty).$$
On the other hand, for each $\alpha\in\R_f$, we have
$P_{\alpha}=(\alpha,0)\in B$ and the corresponding divisor
 $$\div(x-\alpha)=q ((\alpha,0))-q (\infty)=q(P_{\alpha})-q (\infty)$$ is divisible by $q$.
 This implies that
$$a_{P_{\alpha}}= q \cdot b(\alpha)+i.$$
Also, since $\infty$ is neither zero no pole of $h$,
$$0=ni+\sum_{\alpha\in\R_f} b(\alpha) q.$$
Since $n$ and $q$ are relatively prime, $i$ must divide $q$. This
implies that $i=0$ and therefore the divisor
$$D=\div(u(x))=\div(\prod_{\alpha\in \R_f}(x-\alpha)^{b(\alpha)})$$
is divisible by $q$.

Conversely, suppose a divisor $D=\sum_{P\in B}a_P (P)$ with
$\sum_{P\in B}a_P =0$ and all $a_P$ are divisible by $q$. Let us
put
$$h=\prod_{P\in B}(x-x(P))^{a_P/q}.$$
One may easily check that $D=\div(h)$.
\end{proof}

\begin{lem}
\label{cycl}
 $1+\delta_q+ \cdots + \delta_q^{q-1}=0$ in $\End(J(C_{f,q}))$.
The subring $\Z[\delta_q] \subset \End(J(C_{f,q}))$ is isomorphic
to the ring $\Z[t]/{\mathcal P}_q(t)\Z[t]$.  The $\Q$-subalgebra
$\Q[\delta_q]\subset\End^0(J(C_{f,q}))=\End^0(J(C_{f,q}))$ is
isomorphic to $\Q[t]/{\mathcal
P}_q(t)\Q[t]=\prod_{i=1}^r\Q(\zeta_{p^i})$.
\end{lem}

\begin{proof}
If $q=p$ is a prime this assertion is proven in
 \cite[p.~149]{Poonen},  \cite[p.~458]{SPoonen}. So, further we
 may assume that $q>p$. It follows from Remark \ref{nondiv} that we
 may assume that $p$ does {\sl not} divide $n$.

 Now we follow arguments of \cite[p.~458]{SPoonen} (where the case
 of $q=p$ was treated). The group $J(C_{f,q}))(K_a)$ is generated by
 divisor classes of the form $(P)-(\infty)$ where $P$ is a finite
 point on $C_{f,p}$. The divisor of the rational function $x-x(P)$ is
 $(\delta_q^{q-1}P)+\cdots + (\delta_q P) + (P)- q(\infty)$. This
 implies that
 $${\mathcal P}_q(\delta_q)=0 \in \End(J(C_{f,q})).$$
 Applying Remark \ref{alb}(ii), we conclude that ${\mathcal P}_q(t)$
 is the minimal polynomial of $\delta_q$ in $\End(J(C_{f,q}))$.
\end{proof}

Let us define the abelian (sub)variety

$$J^{(f,q)}:= \PP_{q/p}(\delta_q)(J(C_{f,q}))\subset J(C_{f,q}).$$
Clearly, $J^{(f,q)}$ is a $\delta_q$-invariant abelian subvariety
defined over $K(\zeta_q)$. In addition,
$$\Phi_q(\delta_q)(J^{(f,q)})=0.$$

\begin{rem}
\label{qp} If $q=p$ then $\PP_{q/p}(t)=\PP_{1}(t)=1$ and therefore
$J^{(f,p)}=J(C_{f,p})$.
\end{rem}

\begin{rem}
\label{nonzero} Since the polynomials $\Phi_q$ and $\PP_{q/p}$ are
relatively prime, the homomorphism
$$\PP_{q/p}(\delta_q):J^{(f,q)} \to J^{(f,q)}$$
has finite kernel and therefore is an isogeny. In particular, it
is surjective.
\end{rem}

\begin{lem}
\label{fixP}
\begin{enumerate}
\item[(i)] If $p$ does not divide $n$ then
 $\dim (J^{(f,q)})=\frac{(p^r-p^{r-1})(n-1)}{2}$.

 If $q$  divides $n$ then
 $\dim (J^{(f,q)})=\frac{(p^r-p^{r-1})(n-2)}{2}$;
\item[(ii)]
If $p$ does not divide $n$ then
 there is an $K(\zeta_q)$-isogeny $J(C_{f,q})\to
J(C_{f,q/p})\times J^{(f,q)}$.
\item[(iii)]
If $p$ does not divide $n$ and $\zeta \in K$ then the Galois
modules $V_{f,p}$ and
$$(J^{(f,q)})^{\delta_q}:=\{z \in J^{(f,q)}(K_a)\mid \delta_q(z)=z\}$$
are isomorphic.
\end{enumerate}
\end{lem}

\begin{proof}
Clearly, we may assume that $\zeta\in K$.  It follows from
Remark \ref{nondiv} that we
 may assume that $p$ does {\sl not} divide $n$.
Clearly, the assertion (ii) implies the assertion (i). Further we
will prove the assertions (ii) and (iii).

Let us consider the curve
$$C_{f,q/p}:y_1^{q/p}=f(x_1)$$
and a regular surjective map
$$\pi_1:C_{f,q} \to C_{f,q/p}, \quad x_1=x, y_1=y^p.$$
Clearly,
$$\pi_1\delta_q=\delta_{q/p}\pi_1.$$
By Albanese functoriality, $\pi_1$ induces a certain surjective
homomorphism of jacobians $J(C_{f,q}) \twoheadrightarrow
J(C_{f,q/p})$ which we continue to denote by $\pi_1$. Clearly, the
equality $\pi_1\delta_q=\delta_{q/p}\pi_1$ remains true in
$\Hom(J(C_{f,q}),J(C_{f,q/p})$. By Lemma \ref{cycl},
$${\mathcal P}_{q/p}(\delta_{q/p})=0 \in \End(J(C_{f,{q/p}})).$$
It follows from Lemma \ref{nonzero} that
$$\pi_1(J^{(f,q)})=0.$$
It follows that $\dim(J^{(f,q)})$ does not exceed
$$\dim(J(C_{f,q}))-\dim(J(C_{f,q/p}))=
\frac{(p^r-1)(n-1)}{2}-\frac{(p^{r-1}-1)(n-1)}{2}=\frac{(p^r-p^{r-1})(n-1)}{2}.$$

By definition of $J^{(f,q)}$, for each divisor $D=\sum_{P\in B}a_P
(P)$ the linear equivalence class of
$$p^{r-1}D=\sum_{P\in B}p^{r-1}a_P(P)$$
lies in  $(J^{(f,q)})^{\delta_q}\subset J^{(f,q)}(K_a)\subset
J(C_{f,q})(K_a)$. It follows from Lemma \ref{order} that the class
of $p^{r-1}D$ is zero if and only if all $p^{r-1}a_P$ are
divisible by $q=p^r$, i.e. all $a_P$ are divisible by $p$. This
implies that the set of linear equivalence classes of $p^{r-1}D$
is a Galois submodule isomorphic to $V_{f,p}$. We need to prove
that $(J^{(f,q)})^{\delta_q}=V_{f,p}$.

Recall that $J^{(f,q)}$ is $\delta_q$-invariant and the
restriction of $\delta_q$ to $J^{(f,q)}$ satisfies the $q$th
cyclotomic polynomial. This allows us to define the homomorphism
$$\Z[\zeta_q] \to \End(J^{(f,q)})$$
which sends $1$ to the identity map and $\zeta_q$ to $\delta_q$.
Let us put
$$E=\Q(\zeta_q), \O=\Z[\zeta_q]\subset \Q(\zeta_q)=E.$$
It is well-known that $\O$ is the ring of integers in $E$,
$$\lambda=(1-\zeta_q)\Z[\zeta_q]=(1-\zeta_q)\O$$
is a maximal ideal in $\O$ with $\O/\lambda=\F_p$ and
$\O\otimes\Z_p=\Z_p[\zeta_q]$ is the ring of integers in the field
$\Q_p(\zeta_q)$. Notice also that $\O\otimes\Z_p$ coincides with
the completion $\O_{\lambda}$ of $\O$ with respect to
$\lambda$-adic topology and $\O_{\lambda}/\lambda
\O_{\lambda}=\O/\lambda=\F_p$.

 It follows (see \cite{Ribet2}) that
$$d=\frac{2\dim(J^{(f,q)})}{E:\Q]}=\frac{2\dim(J^{(f,q)})}{p^r-p^{r-1}}$$
is a positive integer and the $\Z_p$-Tate module $T_p(J^{(f,q)})$
is a free $\O_{\lambda}$-module of rank $d$. It follows that
$T_p(J^{(f,q)})\otimes_{\O_{\lambda}}\F_p$ is a $d$-dimensional
vector space. On the other hand, clearly
$$(J^{(f,q)})^{\delta_q}=\{u \in J^{(f,q)}(K_a)\mid
(1-\delta_p)(u)=0\}=J^{f,q}_{\lambda}=T_p(J^{f,q})\otimes_{\O_{\lambda}}\F_p.$$
Since $(J^{(f,q)})^{\delta_q}$ contains $(n-1)$-dimensional
$\F_p$-vector space $V_{f,p}$,
$$d \ge n-1.$$
This implies that
$$2\dim(J^{(f,q)})=d (p^r-p^{r-1}) \ge (n-1)(p^r-p^{r-1})$$ and
therefore
$$\dim(J^{(f,q)}) \ge \frac{(n-1)(p^r-p^{r-1})}{2}.$$
But we have already seen that
$$\dim(J^{(f,q)}) \le \frac{(n-1)(p^r-p^{r-1})}{2}.$$
This implies that
$$\dim(J^{(f,q)}) = \frac{(n-1)(p^r-p^{r-1})}{2}.$$
It follows that $d=n-1$ and therefore
$$(J^{(f,q)})^{\delta_q}=V_{f,p}.$$
\end{proof}

\begin{cor}
\label{split} If $p$ does not divide $n$ then there is a
$K(\zeta_q)$-isogeny $J(C_{f,q})\to J(C_{f,p})\times \prod_{i=2}^r
J^{(f,p^i)}=\prod_{i=1}^r J^{(f,p^i)}$.
\end{cor}
\begin{proof}
Combine Corollary \ref{fixP}(ii) and Remark \ref{qp} with easy
induction by $r$.
\end{proof}

\begin{rem}
\label{multJ}
 Suppose that $p$ does not divide $n$ and consider the
induced linear operator
$$\delta_q^*:\Omega^1(J^{(f,q)})\to \Omega^1(J^{(f,q)}).$$
It follows from Theorem \ref{chardif} combined with Remark
\ref{multprim} that its spectrum consists of primitive $q$th roots
of unity $\zeta^{-i}$ ($1\le i<q$) with
$\left[\frac{ni}{q}\right]>0$ and the multiplicity of $\zeta^{-i}$
equals $\left[\frac{ni}{q}\right]$.
\end{rem}

\begin{thm}
\label{maximalV}
 Suppose that $n \ge 5$ is an integer. Let $p$ be a
prime, $r \ge 1$ an integer and $q=p^r$. Suppose  that
 $p$ does not divide $n$. Suppose that  $K$ is a
field of characteristic different from $p$ containing a primitive
$q$th root of unity $\zeta$, Let $f(x)\in K[x]$ be a separable
polynomial of degree $n$ and $\Gal(f)$ its Galois group. Suppose
that the $\Gal(f)$-module $V_{f,p}$ is very simple. Then the image
$\O$ of
$$\Z[\delta_q] \to \End(J^{(f,q)})$$
is isomorphic to $\Z[\zeta_q]$ and
enjoys one of the following two properties.
\begin{itemize}
\item[(i)]
 $\O$
is a maximal commutative subring in $\End(J^{(f,q)})$;
\item[(ii)]
$\fchar(K)>0$ and  the centralizer of $\O\otimes\Q\cong
\Q(\zeta_q)$ in $\End^0(J^{(f,q)})$ is a central simple
$(n-1)^2$-dimensional $\Q(\zeta_q)$-algebra.
\end{itemize}
\end{thm}

\begin{proof}
Clearly, $\O$ is isomorphic to $\Z[\zeta_q]$. Let us put
$\lambda=(1-\zeta_q)\Z[\zeta_q]$. By Lemma \ref{fixP}(iii), the
Galois module $(J^{(f,q)})^{\delta_q}=J^{(f,q)}_{\lambda}$ is
isomorphic to $V_{f,p}$. Applying Theorem \ref{Very}, we conclude
that either (ii) holds true or one of the following conditions
hold:

\begin{itemize}
\item[(a)]
$\O$ is a maximal commutative subring in $\End(J^{(f,q)})$ ;
\item[(b)]
 $\fchar(K)=0$ and there exist a
$\frac{\varphi(q)}{2}$-dimensional abelian variety $Z$ over $K_a$,
an embedding $\Q(\zeta_q) \hookrightarrow \End^0(Z)$ and a
$\Q(\zeta_q)$-equivariant isogeny $\psi:Z^{n-1} \to J^{(f,q)}$.
\end{itemize}

Clearly, if (a) is fulfilled then we are done

If $q=p$ and $\fchar(K)=0$ then it is known \cite{ZarhinTexel},
\cite[Th. ~5.3]{ZarhinCrelle} that  (a) is fulfilled.

So further we may assume  that (b) holds true. In particular,
$\fchar(K)=0$. We may also assume that $q>p$. In order to finish
the proof, we need to arrive to a contradiction. Clearly, $\psi$
induces an isomorphism
$$\psi^*:\Omega^1((J^{(f,q)})) \cong \Omega^1(Z^{n-1})$$
that commutes with the action of $\Q(\zeta_q)$. (Here  again we
use that $\fchar(K)=0$.) Since
$$\dim\Omega^1(Z)=\dim(Z)=\frac{\varphi(q)}{2}.$$
the linear operator in $\Omega^1(Z)$ induced by $\zeta_q$ has, at
most, $\frac{\varphi(q)}{2}$ distinct eigenvalues. It follows that
the linear operator in $\Omega^1(Z^{n-1})=\Omega^1(Z)^{n-1}$
induced by $\zeta_q$ also has, at most, $\frac{\varphi(q)}{2}$
distinct eigenvalues. This implies that the linear operator
$\delta_q^*$ in $\Omega^1((J^{(f,q)}))$ also has, at most,
$\frac{\varphi(q)}{2}$ distinct eigenvalues. Recall that the
eigenvalues of $\delta_q^*$ are primitive $q$th roots of unity
$\zeta^{-i}$ with
$$1\le i<q, (i,p)=1, \left[\frac{ni}{q}\right]>0.$$
Clearly, the inequality $\left[\frac{ni}{q}\right]>0$ means that
$ni \ge q$, i.e.
$$i \ge \frac{q}{n}\ge \frac{q}{5}.$$
So, in order to get a desired contradiction, it suffices that the
cardinality of the set of integers
$$B:=\left\{i\mid \frac{q}{5} \le i < q=p^r, (i,p)=1\right\}$$
is strictly greater than $(p-1)p^{r-1}/2$. Indeed, clearly,
$\frac{p}{5}<\frac{p-1}{2}$ and
$$\#(B)> \varphi(q)-\frac{q}{5}=(p-1)p^{r-1}-\frac{p^{r-1}p}{5}=
(p-1-\frac{p}{5})p^{r-1}>\frac{p-1}{2}p^{r-1}.$$
\end{proof}

\begin{cor}
\label{maximal}
 Suppose that $n \ge 5$ is an integer. Let $p$ be a
prime, $r \ge 1$ an integer and $q=p^r$. Assume in addition that
either $p$ does not divide $n$ or $q\mid n$ and $(n,q)\ne (5,5)$.
Let $K$ be a field of characteristic different from $p$, Let
$f(x)\in K[x]$ be an irreducible separable polynomial of degree
$n$ such that $\Gal(f)=\Sn$ or $\An$. Then the image $\O$ of
$$\Z[\delta_q] \to \End(J^{(f,q)})$$
is isomorphic to $\Z[\zeta_q]$ and
enjoys one of the following two properties.
\begin{itemize}
\item(i)]
 $\O$
is a maximal commutative subring in $\End(J^{(f,q)})$;
\item(ii)]
$\fchar(K)>0$ and the centralizer of $\O\otimes\Q\cong
\Q(\zeta_q)$ in $\End^0(J^{(f,q)})$ is a central simple
$(n-1)^2$-dimensional $\Q(\zeta_q)$-algebra.
\end{itemize}
\end{cor}

\begin{proof}
If $p$ divides $n$ then $n>5$ and therefore $n-1 \ge 5$.  By
Remark \ref{nondiv}, we may assume that $p$ does not divide $n$.
If we replace $K$ by $K(\zeta)$ then  still  $\Gal(f)=\Sn$ or
$\An$. By Remark \ref{SnAn} if $\Gal(f)=\Sn$ or $\An$ then the
$\Gal(f)$-module $V_{f,p}$ is very simple. One has only apply
Theorem \ref{maximalV}.
\end{proof}

\begin{thm}
\label{bigCM} Suppose $n \ge 4$ and $p$ does not divide $n$.
Assume also that $\fchar(K)=0$ and $\Q[\delta_q]$ is a maximal
commutative subalgebra in $\End^0(J^{(f,q)})$. Then
$\End^0(J^{(f,q)})=\Q[\delta_q] \cong\Q(\zeta_q)$ and therefore
 $\End(J^{(f,q)})=\Z[\delta_q] \cong \Z[\zeta_q]$. In particular,
 $J^{(f,q)}$ is an absolutely simple abelian variety.
\end{thm}

\begin{proof} Let  $\CC=\CC_{J^{(f,p)}}$ be the center of
$\End^0(J^{(f,p)})$. Since $\Q[\delta_q]$ is a maximal
commutative, $\CC\subset \Q[\delta_q]$.

Replacing, if necessary, $K$ by  its subfield (finitely) generated
over $\Q$ by all the coefficients of $f$, we may assume that $K$
(and therefore $K_a$) is isomorphic to a subfield of the field
$\C$ of complex numbers. So, $K \subset K_a \subset \C$. We may
also assume that $\zeta=\zeta_q$ and consider $J^{(f,q)}$ as
complex abelian variety.

Let $\Sigma=\Sigma_E$ be the set of all field embeddings
$\sigma:E=\Q[\delta_q]\hookrightarrow \C$. We  are going to apply
Theorem \ref{mult}  to $Z=J^{(f,q)}$ and $E=\Q[\delta_q]$. In
order to do that we need to get some information about the
multiplicities
$$n_{\sigma}=\n_{\sigma}(Z,E)=n_{\sigma}(J^{(f,q)},\Q[\delta_q]).$$
Remark \ref{dual} allows us to do it, using the action of
$\Q[\delta_q]$  on the space $\Omega^1(J^{(f,q)})$ of
differentials of the first kind on $J^{(f,q)}$.

In other words, $\Omega^1(J^{(f,q)})_{\sigma}$ is the eigenspace
corresponding to the eigenvalue $\sigma(\delta_q)$ of $\delta_q$
and $n_{\sigma}$ is the multiplicity of the eigenvalue
$\sigma(\delta_q)$.

 Let $i<q$ be a positive integer that is not divisible by $p$ and
$\sigma_i:\Q[\delta_p]\hookrightarrow \C$ be the embedding which
sends $\delta_p$ to $\zeta^{-i}$.  Clearly, for each $\sigma$
there exists precisely one $i$ such that $\sigma=\sigma_i$.
Clearly, $\Omega^1(J^{(f,q)})_{\sigma_i}$ is the eigenspace of
$\Omega^1(J^{(f,q)})$ attached to the eigenvalue $\zeta^{-i}$ of
$\delta_q$. Therefore $n_{\sigma_i}$ coincides with the
multiplicity of the eigenvalue $\zeta^{-i}$. It follows from
Remark \ref{multJ} that
$$n_{\sigma_i}=\left[\frac{ni}{q}\right].$$
Now the assertion of the Theorem follows from Corollary
\ref{cyclmult} applied to $E=\Q(\zeta_q)\cong\Q[\delta_q]$.
\end{proof}

Combining Corollary \ref{maximal} and \ref{bigCM}, we obtain the
following statement.

\begin{thm}
\label{final} Let $p$ be a prime, r a positive integer, $q=p^r$
and $K$ a field of characteristic zero. Suppose that $f(x)\in
K[x]$ is an irreducible polynomial of degree $n\ge 5$ and
$\Gal(f)=\Sn$ or $\An$. Assume also that either $p$ does not
divide $n$ or $q$ divides $n$. Then
 $\End^0(J^{(f,q)})=\Q[\delta_q]\cong \Q(\zeta_q)$ and therefore
$\End(J^{(f,q)})=\Z[\delta_q]\cong \Z[\zeta_q]$. In particular,
 $J^{(f,q)}$ is an absolutely simple abelian variety.
\end{thm}

Combining Theorem \ref{bigCM} amd Corollary \ref{maximalV}, we
obtain the following statement.

\begin{thm}
\label{handysup} Let $p$ be a prime, r a positive integer, $q=p^r$
and $K$ a field of characteristic zero. Let $f(x)\in K[x]$ be an
polynomial of degree $n\ge 5$. Assume also that $p$ does not
divide $n$ and the $\Gal(f)$-module $V_{f,p}$ is very simple. Then
 $\End^0(J^{(f,q)})=\Q[\delta_q]\cong \Q(\zeta_q)$ and therefore
$\End(J^{(f,q)})=\Z[\delta_q]\cong \Z[\zeta_q]$. In particular,
 $J^{(f,q)}$ is an absolutely simple abelian variety.
\end{thm}

\section{Jacobians and their endomorphism rings}
\label{prf}
 Throughout this section we assume that $K$ is a field of characteristic zero. Recall that $K_a$ is
an algebraic closure of $K$ and $\zeta \in K_a$ is a primitive
$q$th root of unity. Suppose  $f(x)\in K[x]$ is a polynomial of
degree $n \ge 5$ without multiple roots, $\R_f \subset K_a$ is the
set of its roots, $K(\R_f)$ is its splitting field. Let us put
$$\Gal(f)=\Gal(K(\R_f)/K)\subset\Perm(\R_f).$$
Let $r$ be a positive integer. Recall (Corollary \ref{split}) that
if  $p$ does not divide $n$ then there is a
$K(\zeta_{p^r})$-isogeny $J(C_{f,p^r})\to  \prod_{i=1}^r
J^{(f,p^i)}$. Applying Theorem \ref{handysup} to all $q=p^i$, we
obtain the following assertions.

\begin{thm}
\label{handysupJ} Let $p$ be a prime, r a positive integer,
$q=p^r$ and $K$ a field of characteristic zero. Let $f(x)\in K[x]$
be an polynomial of degree $n\ge 5$. Assume also that $p$ does not
divide $n$ and the $\Gal(f)$-module $V_{f,p}$ is very simple. Then
$$\End^0(J(C_{f,q}))=\Q[\delta_q]
\cong\Q[t]/{\mathcal P}_q(t)\Q[t]=\prod_{i=1}^r\Q(\zeta_{p^i}).$$
\end{thm}

\begin{thm}
\label{finalJ} Let $p$ be a prime, r a positive integer and $K$ a
field of characteristic zero. Suppose that $f(x)\in K[x]$ is an
irreducible polynomial of degree $n\ge 5$ and $\Gal(f)=\Sn$ or
$\An$. Assume also that  either $p$ does not divide $n$ or
$$\End^0(J(C_{f,q}))=\Q[\delta_q]
\cong\Q[t]/{\mathcal P}_q(t)\Q[t]=\prod_{i=1}^r\Q(\zeta_{p^i}).$$
\end{thm}

\begin{proof}
The existence of the isogeny $J(C_{f,q})\to \prod_{i=1}^r
J^{f,(p^i)}$ combined with Theorem \ref{final} implies that the
assertion holds true if $p$ does not divide $n$. If  $q$ divides
$n$ then  Remark \ref{nondiv} allows us to reduce this case to the
already proven case when $p$ does not divide $n-1$.
\end{proof}

\begin{ex}
Suppose $L=\C(z_1, \cdots , z_n)$ is the field of rational
functions in $n$ independent variables $z_1, \cdots , z_n$ with
constant field $\C$ and $K=L^{\Sn}$ is the subfield of symmetric
functions. Then $K_a=L_a$ and $$f(x)=\prod_{i=1}^n(x-z_i) \in
K[x]$$ is an irreducible polynomial over $K$ with Galois group
$\Sn$. Let Let $q=p^r$ be a power of a prime $p$. Let $C$ be a
smooth projective model of the $K$-curve $y^q=f(x)$ and $J(C)$ its
jacobian. It follows from  Theorem \ref{finalJ} that if $n\ge 5$
and either $p$ does not divide $n$ or $q$ divides $n$ then the
algebra of $L_a$-endomorphisms of $J(C)$ is
$\prod_{i=1}^r\Q(\zeta_{p^i})$.
\end{ex}

\begin{ex}
Let $h(x) \in \C[x]$ be a  {\sl Morse polynomial} of degree $n \ge
5$. This means that the derivative $h'(x)$ of $h(x)$ has $n-1$
distinct roots $\beta_1, \cdots \beta_{n-1}$  and $h(\beta_i) \ne
h(\beta_j)$ while $i\ne j$. (For example, $x^n-x$ is a Morse
polynomial.) Let $K=\C(z)$ be the field of rational functions in
variable $z$ with constant field $\C$ and $K_a$ its algebraic
closure. Then a theorem of Hilbert  (\cite[theorem~4.4.5,
p.~41]{Serre}) asserts that the Galois group of $h(x)-z$ over
$k(z)$ is $\Sn$. Let $q=p^r$ be a power of a prime $p$. Let $C$ be
a smooth projective model of the $K$-curve $y^q=h(x)-z$ and $J(C)$
its jacobian. It follows from Theorem \ref{finalJ} that if either
$p$ does not divide $n$ or $q$ divides $n$ then the algebra of
$K_a$-endomorphisms of $J(C)$ is $\prod_{i=1}^r\Q(\zeta_{p^i})$.
\end{ex}

\end{document}